\documentclass{amsproc}%
\usepackage{graphicx}
\usepackage{amscd}
\usepackage{amsmath}
\usepackage{amsfonts}
\usepackage{amssymb}%
\theoremstyle{plain}

\numberwithin{equation}{section}

\begin{document}
\title[Power-conjugate equations in symmetric groups]{Power-conjugate equations in symmetric groups}
\author{Szilvia Homolya}
\address{Institute of Mathematics, University of Miskolc, \\
3515 Miskolc-Egyetemv\'{a}ros, Hungary}
\email{szilvia.homolya@uni-miskolc.hu}
\author{Jen\H{o} Szigeti}
\address{Institute of Mathematics, University of Miskolc, \\
3515 Miskolc-Egyetemv\'{a}ros, Hungary}
\email{matjeno@uni-miskolc.hu}
\thanks{The second named author was partially supported by the National Research,
Development and Innovation Office of Hungary (NKFIH) K138828.}

\begin{abstract}
First we consider the solutions of the general "cubic" equation%
\[
\alpha_{1}\circ x^{r_{1}}\circ\alpha_{2}\circ x^{r_{2}}\circ\alpha_{3}\circ
x^{r_{3}}=1
\]

\noindent(with $r_{1},r_{2},r_{3}\in\{1,-1\}$) in the symmetric group
$\mathrm{S}_{n}$. In certain cases this equation can be rewritten as
$\alpha\circ y\circ\alpha^{-1}=y^{2}$ or as $\alpha\circ y\circ\alpha
^{-1}=y^{-2}$, where $\alpha\in\mathrm{S}_{n}$ depends on the $\alpha_{i}$'s
and the new unknown permutation $y\in\mathrm{S}_{n}$ is a product of $x$ (or
$x^{-1}$) and one of the permutations $\alpha_{i}^{\pm1}$. Using combinatorial
arguments and some basic number theoretical facts, we obtain results about the
solutions of the so-called power-conjugate equation $\alpha\circ y\circ
\alpha^{-1}=y^{e}$ in $\mathrm{S}_{n}$, where $e\in\mathbb{Z}$ is an integer
exponent. A divisibility condition involving the type of $\alpha$ provide
solutions with $\alpha\circ y\neq y\circ\alpha$ and a further condition gives
the complete list of solutions. Some other divisibility assumptions concerning
the type of $\alpha$ ensure that the solutions of $\alpha\circ y\circ
\alpha^{-1}=y^{e}$\ are exactly the solutions of $y^{e-1}=1$ in the
centralizer of $\alpha$. Slightly stronger assumptions provide a complete
answer to the question, when our equation has only the trivial solution $y=1$.

\end{abstract}
\maketitle

\noindent1. INTRODUCTION

\bigskip

\noindent One of the starting points of classical algebra is the solution of
polynomial equations in fields. Thus, the investigation of equations in groups
is a natural idea. Recently a good number of publications appeared that are
related to complexity, taking an algorithmic approach to solving such
equations. In the present paper we are interested only in the explicit
solutions, so we restrict our consideration to the non-algorithmic aspects. It
is a surprising fact that the authors found only a limited number of results
about the explicit solutions of equations in groups. One of the earliest
results (due to Frobenius) is about the number of solutions of $x^{m}=1$ in a
finite group (see [IR]). Further results concerning the equation $x^{m}=1$ in
the symmetric group $\mathrm{S}_{n}$ consisting of all bijective
$\{1,2,\ldots,n\}\longrightarrow\{1,2,\ldots,n\}$ functions can be found in
[CHS], [MW] and [FM]. A general equation (containing constants and group
operations) for a single unknown permutation $x\in\mathrm{S}_{n}$ is of the
form:%
\[
\alpha_{1}\circ x^{r_{1}}\circ\alpha_{2}\circ\cdots\circ\alpha_{k}\circ
x^{r_{k}}=1,
\]
where $k\geq1$, $\alpha_{i}\in\mathrm{S}_{n}$ and $r_{i}\in\{1,-1\}$ for each
$1\leq i\leq k$. Our first impression is that the complete solution of the
above equation is hopeless, on the other hand to deal with some special cases
seems to be a challenging problem. We note that in the above equation
$r_{1}=1$ can be assumed, otherwise $\alpha_{1}\circ x^{-r_{1}}\circ\alpha
_{2}\circ\cdots\circ\alpha_{k}\circ x^{-r_{k}}=1$ is the same equation for the
inverse $x^{-1}$ and $-r_{1}=1$.

\noindent The solutions in the "quadratic" case $k=2$ can easily be obtained
by a simple procedure. If $1=r_{1}=r_{2}$, then $\alpha_{1}\circ x\circ
\alpha_{2}\circ x=1$ is equivalent to%
\[
(x\circ\alpha_{2})^{2}=\alpha_{1}^{-1}\circ\alpha_{2}.
\]
The above square root equation has a solution if and only if the number of the
$2i$-element cycles in $\alpha_{1}^{-1}\circ\alpha_{2}$\ is even for all
integers $i\geq1$. If $1=r_{1}=-r_{2}$, then $\alpha_{1}\circ x\circ\alpha
_{2}\circ x^{-1}=1$ is equivalent to%
\[
x\circ\alpha_{2}\circ x^{-1}=\alpha_{1}^{-1}.
\]
The above equation has a solution if and only if the permutations $\alpha_{2}$
and $\alpha_{1}^{-1}$\ are of the same type. The type of a permutation $\pi
\in\mathrm{S}_{n}$\ is a sequence $\mathrm{type}(\pi)=\left\langle t_{1}%
,t_{2},\ldots,t_{n}\right\rangle $ of integers, where $t_{i}\geq0$ denotes the
number of cycles in $\pi$ of length $i\geq1$. A well-known fact is that for
$\pi_{1},\pi_{2}\in\mathrm{S}_{n}$\ the equality $\mathrm{type}(\pi
_{1})=\mathrm{type}(\pi_{2})$ is equivalent to the conjugate relation between
$\pi_{1}$ and $\pi_{2}$ (there exists a permutation $\tau\in\mathrm{S}_{n}$
such that $\tau\circ\pi_{1}\circ\tau^{-1}=\pi_{2}$). Further details (about
the complete solutions in the case $k=2$) are left to the readers.

\noindent The situation in the case $k\geq3$ is far more complicated. A nice
summary about the general situation can be found in [L]. An important
direction of research is to find solutions of a given group-equation in an
appropriate extension of the base group. The Kervaire--Laudenbach (KL)
conjecture asserts that if the length $r_{1}+r_{2}+\cdots+r_{k}$ of the
equation $\alpha_{1}\circ x^{r_{1}}\circ\alpha_{2}\circ\cdots\circ\alpha
_{k}\circ x^{r_{k}}=1$ over an arbitrary group $G$\ is nonzero, then this
equation has a solution in a group $H$ containing $G$ (here we assume that
$r_{i-1}+r_{i}\neq0$ if $\alpha_{i}=1$, $2\leq i\leq k$). A consequence of a
general extension theorem of Gerstenhaber and Rothaus (see [GR]) is that (KL)
holds for finite groups. For $k=5$ the conjecture (KL) is proved in [E].

\noindent Now consider the general "cubic" equation%
\[
(\ast)\text{ \ \ \ \ \ \ \ \ \ \ \ }\alpha_{1}\circ x^{r_{1}}\circ\alpha
_{2}\circ x^{r_{2}}\circ\alpha_{3}\circ x^{r_{3}}=1
\]
in $\mathrm{S}_{n}$, where $r_{1}=1$. According to the choice of the exponents
we have the following four possibilities%
\[
(\ast1)\text{ \ }\alpha_{1}\circ x\circ\alpha_{2}\circ x\circ\alpha_{3}\circ
x^{-1}=1,\text{ \ }(\ast2)\text{ \ }\alpha_{1}\circ x\circ\alpha_{2}\circ
x^{-1}\circ\alpha_{3}\circ x=1,
\]%
\[
(\ast3)\ \alpha_{1}\circ x\circ\alpha_{2}\circ x^{-1}\circ\alpha_{3}\circ
x^{-1}=1,\ (\ast4)\ \alpha_{1}\circ x\circ\alpha_{2}\circ x\circ\alpha
_{3}\circ x=1.
\]
Clearly, the above equations can be rewritten as follows%
\[
(\ast1)\text{ }\alpha\circ y\circ\beta=y^{2}\text{, where }y=x\circ\alpha
_{2}\text{ and }\alpha=\alpha_{1}^{-1},\beta=\alpha_{2}\circ\alpha_{3}%
^{-1}\circ\alpha_{2}^{-1},
\]%
\[
(\ast2)\text{ }\alpha\circ y\circ\beta=y^{2}\text{, where }y=x^{-1}\circ
\alpha_{1}^{-1}\text{ and }\alpha=\alpha_{2},\beta=\alpha_{1}\circ\alpha
_{3}\circ\alpha_{1}^{-1},
\]%
\[
(\ast3)\text{ }\alpha\circ y\circ\beta=y^{2}\text{, where }y=x\circ\alpha
_{3}^{-1}\text{ and }\alpha=\alpha_{1},\beta=\alpha_{3}\circ\alpha_{2}%
\circ\alpha_{3}^{-1},
\]%
\[
(\ast4)\text{ }\alpha\circ y\circ\beta=y^{-2}\text{, where }y=\alpha_{3}\circ
x\text{ and }\alpha=\alpha_{1}\circ\alpha_{3}^{-1},\beta=\alpha_{2}\circ
\alpha_{3}^{-1}.
\]
Thus, the solution of the "cubic" equation can be reduced to the solution of
the equations%
\[
\alpha\circ y\circ\beta=y^{2}\text{ and }\alpha\circ y\circ\beta
=y^{-2}\text{,}%
\]
where $\alpha,\beta\in\mathrm{S}_{n}$ are constants and the new unknown
permutation is $y\in\mathrm{S}_{n}$. In both cases we assume that $y=1$\ is a
solution of the given equation. Our requirement is quite natural, however the
weaker assumption, that we have at least one solution in $\mathrm{S}_{n}$ is
even more natural. The condition that $y=1$\ is a solution is equivalent to
the fact that $\beta=\alpha^{-1}$ is the inverse of $\alpha$ (in both cases).
In view of the above observations, we consider the power-conjugate equation%
\[
(e\ast\alpha):\text{ }\alpha\circ y\circ\alpha^{-1}=y^{e},
\]
where $e\in\mathbb{Z}$ is an integer exponent. If $e\in\{-1,1\}$, then
$(e\ast\alpha)$ is quadratic and the case $e=0$ is trivial. Therefore in the
rest of the paper we assume that $e\notin\{-1,0,1\}$. Using combinatorial
arguments and some basic number theoretical facts, in certain cases (depending
on the type of $\alpha$) we are able to obtain essential information about the
solutions of $(e\ast\alpha)$\ in $\mathrm{S}_{n}$. Solutions with $\alpha\circ
y\neq y\circ\alpha$\ appear in Theorems 3.1, 3.6 as well as in Corollaries
3.2, 3.3, while solutions with $\alpha\circ y=y\circ\alpha$\ appear in
Theorems 4.5 and 4.8. One of our main results reveals that under certain mild
divisibility assumptions involving the type of $\alpha$, the solutions of
$(e\ast\alpha)$ are exactly the solutions of $y^{e-1}=1$ in the centralizer of
$\alpha$ (see 4.5). Using slightly stronger assumptions, we obtain a complete
answer to the question, when $(e\ast\alpha)$ has only the trivial solution
$y=1$ (see 4.8).

\bigskip

\noindent2.\thinspace PRELIMINARY LEMMAS

\bigskip

\noindent Let $\delta=\tau\circ\alpha\circ\tau^{-1}$ be a conjugate of
$\alpha$ (here $\tau\in\mathrm{S}_{n}$ is fixed). Since the conjugation is an
automorphism of $\mathrm{S}_{n}$, the solutions of $(e\ast\delta)$:
$\delta\circ z\circ\delta^{-1}=z^{e}$ can be obtained as the conjugates
$z=\tau\circ y\circ\tau^{-1}$ of the solutions of $(e\ast\alpha)$. Thus, the
number of solutions of $(e\ast\alpha)$ depends only on the conjugacy class (or
the type) of $\alpha$.

\noindent For a given permutation $\alpha\in\mathrm{S}_{n}$ with
$\mathrm{type}(\alpha)=\left\langle g_{1},g_{2},\ldots,g_{n}\right\rangle
$\ and for an integer $1\leq d\leq n$ we define a set of integers (the
$d$-range of $\alpha$) as%
\[
F_{d}(\alpha)=\left\{  \left.  \underset{1\leq j\leq n,\;d\mid j}{\sum}%
q_{j}\cdot j\right\vert 0\leq q_{j}\leq g_{j}\text{ for all }1\leq j\leq
n\right\}  =
\]%
\[
\{q_{d}\cdot d+q_{2d}\cdot2d+\cdots+q_{\left\lfloor n/d\right\rfloor d}%
\cdot\left\lfloor n/d\right\rfloor d\mid0\leq q_{id}\leq g_{id}\text{ for all
}1\leq i\leq\left\lfloor n/d\right\rfloor \}.
\]
Now $g_{1}\cdot1+g_{2}\cdot2+\cdots+g_{n}\cdot n=n$ gives that $F_{d}%
(\alpha)\subseteq\{d,2d,\ldots,\left\lfloor n/d\right\rfloor d\}$ and
$g_{n}=1$\ implies that $g_{1}=g_{2}=\cdots=g_{n-1}=0$, whence $F_{1}%
(\alpha)=\{0,n\}$\ follows. Clearly, the divisibility $d_{1}\mid d_{2}%
$\ implies that $F_{d_{2}}(\alpha)\subseteq F_{d_{1}}(\alpha)$.

\bigskip

\noindent\textbf{2.1. Lemma.}\textit{ If the containment }$\alpha(H)\subseteq
H$\textit{ holds for some subset }$H\subseteq\{1,2,\ldots,n\}$\textit{, then
}$\alpha(H)=H$\textit{ and }$H=D_{1}\cup D_{2}\cup\cdots\cup D_{m}$\textit{ is
a union of certain pairwise disjoint cycles of }$\alpha$\textit{. If }$1\leq
d\leq n$\textit{ and }$d\mid\left\vert D_{j}\right\vert $\textit{ holds for
all }$1\leq j\leq m$\textit{, then for the number of elements we have}%
\[
\left\vert H\right\vert =\left\vert D_{1}\right\vert +\left\vert
D_{2}\right\vert +\cdots+\left\vert D_{m}\right\vert \in F_{d}(\alpha).
\]
\noindent\textbf{Proof.} Obvious. $\square$

\bigskip

\noindent\textbf{2.2. Lemma.}\textit{ If }$\alpha\circ y\circ\alpha^{-1}%
=z$\textit{ holds for }$\alpha,y,z\in\mathrm{S}_{n}$\textit{ and }%
$(c_{1},c_{2},\ldots,c_{r})$\textit{ is a cycle of }$y$\textit{, then
}$(\alpha(c_{1}),\alpha(c_{2}),\ldots,\alpha(c_{r}))$\textit{ is a cycle of
}$z$\textit{.}

\bigskip

\noindent\textbf{Proof.} Clearly, $z(\alpha(c_{i}))=(z\circ\alpha
)(c_{i})=(\alpha\circ y)(c_{i})=\alpha(y(c_{i}))=\alpha(c_{i+1})$ (indices
$1\leq i\leq r$ are considered modulo $r$). $\square$

\bigskip

\noindent\textbf{2.3. Lemma.}\textit{ Consider the union }$\{1,2,\ldots
,n\}=H_{1}\cup H_{2}\cup\cdots\cup H_{s}$\textit{ of the pairwise disjoint
fixed subsets }$H_{k}\subseteq\{1,2,\ldots,n\}$\textit{, }$1\leq k\leq
s$\textit{ of }$\alpha\in\mathrm{S}_{n}$\textit{ (now }$\alpha(H_{k})=H_{k}%
$\textit{\ for each }$1\leq k\leq s$\textit{)\ and let }$y_{k}\in
\mathrm{S}_{H_{k}}$\textit{ (here }$y_{k}:H_{k}\longrightarrow H_{k}$\textit{
is a permutation) be a solution of }$(e\ast(\alpha\upharpoonright H_{k}%
))$\textit{, where }$(\alpha\upharpoonright H_{k}):H_{k}\longrightarrow H_{k}%
$\textit{ is the restriction of }$\alpha$\textit{\ to }$H_{k}$\textit{. Now
the disjoint union }$y_{1}\sqcup y_{2}\sqcup\cdots\sqcup y_{s}$\textit{ of the
permutations }$y_{k}$\textit{ (}$1\leq k\leq s$\textit{) is a solution of
}$(e\ast\alpha)$\textit{\ in }$\mathrm{S}_{n}$\textit{. Notice that for }%
$i\in\{1,2,\ldots,n\}$%
\[
(y_{1}\sqcup y_{2}\sqcup\cdots\sqcup y_{s})(i)=y_{k}(i)\text{\textit{, where
}}1\leq k\leq s\text{\textit{ is the unique index with }}i\in H_{k}%
\text{\textit{.}}%
\]

\bigskip

\noindent\textbf{Proof.} Obvious. $\square$

\bigskip

\noindent\textbf{2.4. Lemma.}\textit{ If }$y_{1},y_{2}\in\mathrm{S}_{n}%
$\textit{ are solutions of }$(e\ast\alpha)$\textit{ such that }$y_{1}\circ
y_{2}=y_{2}\circ y_{1}$\textit{, then the product }$y_{1}\circ y_{2}$\textit{
is also a solution of }$(e\ast\alpha)$.\textit{ If }$y\in\mathrm{S}_{n}%
$\textit{ is a solution of }$(e\ast\alpha)$\textit{, then }$y^{-1}%
\in\mathrm{S}_{n}$\textit{ is also a solution of }$(e\ast\alpha)$\textit{ and
for all }$i\in\mathbb{Z}$\textit{ and }$k\geq0$\textit{ we have}%
\[
\alpha^{k}\circ y^{i}=y^{e^{k}i}\circ\alpha^{k}.
\]

\bigskip

\noindent\textbf{Proof.} Since $\alpha\circ y_{1}\circ\alpha^{-1}=y_{1}^{e}$
and $\alpha\circ y_{2}\circ\alpha^{-1}=y_{2}^{e}$ hold, the multiplicative
property of the conjugation gives that%
\[
\alpha\circ y_{1}\circ y_{2}\circ\alpha^{-1}=(\alpha\circ y_{1}\circ
\alpha^{-1})\circ(\alpha\circ y_{2}\circ\alpha^{-1})=y_{1}^{e}\circ y_{2}%
^{e}=(y_{1}\circ y_{2})^{e}.
\]
If $y\in\mathrm{S}_{n}$ is a solution of $(e\ast\alpha)$, then%
\[
\alpha\circ y^{-1}\circ\alpha^{-1}=(\alpha\circ y\circ\alpha^{-1})^{-1}%
=(y^{e})^{-1}=(y^{-1})^{e}.
\]
Since $y^{i}$ is also a solution of $(e\ast\alpha)$, we have $\alpha\circ
y^{i}\circ\alpha^{-1}=(y^{i})^{e}$, whence $\alpha\circ y^{i}=y^{ei}%
\circ\alpha$ follows for all $i\geq0$. For $k\geq0$ we use an induction:%
\[
\alpha^{k+1}\circ y^{i}\!=\!\alpha\circ(\alpha^{k}\circ y^{i})\!=\!\alpha
\circ(y^{e^{k}i}\circ\alpha^{k})\!=\!(\alpha\circ y^{e^{k}i})\circ\alpha
^{k}\!=\!(y^{e(e^{k}i)}\circ\alpha)\circ\alpha^{k}\!=\!y^{e^{k+1}i}\circ
\alpha^{k+1}\!.\!\square
\]

\bigskip

\noindent\textbf{2.5. Lemma.}\textit{ If }$y\in\mathrm{S}_{n}$\textit{ is a
solution of }$(e\ast\alpha)$\textit{, then }$y$\textit{ and }$y^{e}$\textit{
are of the same type }$\left\langle t_{1},t_{2},\ldots,t_{n}\right\rangle
$\textit{, }$y^{e^{w}-1}=1$\textit{ and for any cycle }$(c_{1},c_{2}%
,\ldots,c_{r})$\textit{ of }$y$\textit{ the cycle length }$r\geq1$\textit{ is
a divisor of }$e^{w}-1$\textit{ (i.e. }$r\mid e^{w}-1$\textit{), where
}$\mathrm{type}(\alpha)=\left\langle g_{1},g_{2},\ldots,g_{n}\right\rangle
$\textit{ and}%
\[
w=\mathrm{ord}(\alpha)=\operatorname{lcm}\{a\mid1\leq a\leq n,g_{a}\neq0\}
\]
\textit{is the order of }$\alpha$\textit{. We also have }$\gcd(r,e)=1$%
\textit{\ and }$(c_{1},c_{e+1},c_{2e+1},\ldots,c_{(r-1)e+1})$\textit{ is a
cycle of }$y^{e}$\textit{, where the indices }$ie+1$\textit{, }$0\leq i\leq
r-1$\textit{ are taken in }$\{1,2,\ldots,r\}$\textit{ modulo }$r$\textit{.
This cycle of }$y^{e}$\textit{ has the same elements as the original
cycle\ and each cycle of }$y^{e}$\textit{\ can be obtained by the above
construction, starting from a uniquely determined cycle of }$y$\textit{.}

\bigskip

\noindent\textbf{Proof.} According to $(e\ast\alpha)$, the permutation $y^{e}$
is the conjugate of $y$ (by $\alpha$), whence $\mathrm{type}(y)=\mathrm{type}%
(y^{e})$ follows. The application of Lemma 2.4 gives that $\alpha^{w}\circ
y=y^{e^{w}}\circ\alpha^{w}$. Now $\alpha^{w}=1$ implies that $y^{e^{w}-1}=1$,
whence we obtain that any cycle length $r\geq1$\ of $y$\ is a divisor of
$e^{w}-1$. Clearly, $\gcd(r,e)=1$\ is a consequence of $r\mid e^{w}-1$. If
$(c_{1},c_{2},\ldots,c_{r})$ is a cycle of $y$, then $(c_{1},c_{e+1}%
,c_{2e+1},\ldots,c_{(r-1)e+1})$ is obviously a cycle of $y^{e}$ of the same
length $r$ (notice that $c_{ie+1}\neq c_{je+1}$ follows from $r\nmid(j-i)e$
for all $1\leq i<j\leq r-1$). If $(c_{1},y^{e}(c_{1}),\ldots,y^{(r-1)e}%
(c_{1}))$ is a cycle of $y^{e}$ of length $r\geq1$, then $y^{re}(c_{1})=c_{1}$
and $s\mid re$, where $y^{s}(c_{1})=c_{1}$ and $1\leq s\leq n$ is the length
of the $y$-cycle strating with $c_{1}$. Since $\gcd(s,e)=1$, we obtain that
$s\mid r$. The containment $\{c_{1},y^{e}(c_{1}),\ldots,y^{(r-1)e}%
(c_{1})\}\subseteq\{c_{1},y(c_{1}),\ldots,y^{s-1}(c_{1})\}$\ implies that
$r\leq s$, whence $r=s$ follows. Thus, the above cycle of $y^{e}$\ can be
constructed by the given process starting from the cycle $(c_{1}%
,y(c_{1}),\ldots,y^{r-1}(c_{1}))$ of $y$. $\square$

\bigskip

\noindent\textbf{2.6. Lemma.}\textit{ Let }$y\in\mathrm{S}_{n}$\textit{ be a
solution of }$(e\ast\alpha)$\textit{ with }$\mathrm{type}(y)=\left\langle
t_{1},t_{2},\ldots,t_{n}\right\rangle $\textit{ and for a given }$1\leq r\leq
n$\textit{ with }$t_{r}\geq1$\textit{ consider the base sets}%
\[
C_{i}^{(r)}\subseteq\{1,2,\ldots,n\},1\leq i\leq t_{r}%
\]
\textit{of all }$r$\textit{-element cycles of }$y$\textit{ (}$\left\vert
C_{i}^{(r)}\right\vert =r$\textit{ for all }$1\leq i\leq t_{r}$\textit{). Then
there exists a permutation }$\gamma$\textit{\ of the index set }%
$\{1,2,\ldots,t_{r}\}$\textit{ such that}%
\[
\alpha(C_{i}^{(r)})=C_{\gamma(i)}^{(r)}.
\]
\textit{Since }$\gamma$\textit{\ is uniquely determined by }$\alpha$\textit{,
it is natural to use the notation }$\gamma=\alpha^{(r)}$\textit{. Now }%
$t_{r}r\in F_{1}(\alpha)$\textit{\ and if }$\alpha^{(r)}$\textit{\ has a cycle
}$(i_{1},i_{2},\ldots,i_{d})$\textit{\ of length }$1\leq d\leq t_{r}$\textit{,
then }$d\mid w=\mathrm{ord}(\alpha)$\textit{ and }$dr\in F_{d}(\alpha
)$\textit{.}

\bigskip

\noindent\textbf{Proof.} In view of Lemmas 2.2 and 2.5, any $r$-element cycle
of $y^{e}$ can be obtained as the $\alpha$ image of a unique $r$-element cycle
of $y$. Since the base sets of the $r$-element cycles in $y$ and in $y^{e}$
coincide, we obtain that $\alpha(C_{i}^{(r)})=C_{\gamma(i)}^{(r)}$ for some
permutation $\gamma$\ of the indices. Now $H=C_{1}^{(r)}\cup C_{2}^{(r)}%
\cup\cdots\cup C_{t_{r}}^{(r)}$ is a fixed set of $\alpha$, whence
$t_{r}r=\left\vert H\right\vert \in F_{1}(\alpha)$ follows by Lemma 2.1. If
$(i_{1},i_{2},\ldots,i_{d})$ is a cycle of $\alpha^{(r)}$, then%
\[
\alpha(C_{i_{1}}^{(r)})=C_{i_{2}}^{(r)},\alpha(C_{i_{2}}^{(r)})=C_{i_{3}%
}^{(r)},\ldots,\alpha(C_{i_{d-1}}^{(r)})=C_{i_{d}}^{(r)},\alpha(C_{i_{d}%
}^{(r)})=C_{i_{1}}^{(r)}.
\]
Now $H(d)=C_{i_{1}}^{(r)}\cup C_{i_{2}}^{(r)}\cup\cdots\cup C_{i_{d}}^{(r)}$
is a fixed set of $\alpha$ and the above property of $\alpha$ (or
$\alpha^{(r)})$ gives that%
\[
H(d)=D_{1}\cup D_{2}\cup\cdots\cup D_{m},
\]
where $D_{1},D_{2},\ldots,D_{m}$ are certain pairwise disjoint cycles of
$\alpha$ such that $d\mid\left\vert D_{j}\right\vert $ and $\left\vert
D_{j}\right\vert \mid w$\ hold for all $1\leq j\leq m$. Thus, $d\mid w$ and%
\[
dr=\left\vert H(d)\right\vert =\left\vert D_{1}\right\vert +\left\vert
D_{2}\right\vert +\cdots+\left\vert D_{m}\right\vert \in F_{d}(\alpha)
\]
follows by the repeated application of Lemma 2.1. $\square$

\bigskip

\noindent\textbf{2.7. Lemma.}\textit{ Let }$y\in\mathrm{S}_{n}$\textit{ be a
solution of }$(e\ast\alpha)$\textit{ with }$\mathrm{type}(y)=\left\langle
t_{1},t_{2},\ldots,t_{n}\right\rangle $\textit{. If }$t_{r}\neq0$\textit{ and
the induced permutation }$\alpha^{(r)}$\textit{\ of the indices }%
$\{1,2,\ldots,t_{r}\}$\textit{ (see Lemma 2.6) has a cycle }$(i_{1}%
,i_{2},\ldots,i_{d})$\textit{\ of length }$1\leq d\leq t_{r}$\textit{ and
}$\gcd(e^{d}-1,r)=1$\textit{, then }$\alpha$\textit{\ also has a cycle of
length }$d$\textit{ in }$C_{i_{1}}^{(r)}\cup C_{i_{2}}^{(r)}\cup\cdots\cup
C_{i_{d}}^{(r)}$\textit{.}

\bigskip

\noindent\textbf{Proof.} Let $C_{i_{1}}^{(r)}=\{c_{1},c_{2},\ldots,c_{r}\}$ be
the set of all elements in a cycle $(c_{1},c_{2},\ldots,c_{r})$ of $y$ (notice
that $r\mid e^{w}-1$ by Lemma 2.5). Since $c_{1}\in C_{i_{1}}^{(r)}$ and%
\[
\alpha(C_{i_{1}}^{(r)})=C_{i_{2}}^{(r)},\alpha(C_{i_{2}}^{(r)})=C_{i_{3}%
}^{(r)},\ldots,\alpha(C_{i_{d-1}}^{(r)})=C_{i_{d}}^{(r)},\alpha(C_{i_{d}%
}^{(r)})=C_{i_{1}}^{(r)},
\]
the containment $\alpha^{d}(c_{1})\in C_{i_{1}}^{(r)}$ holds, whence
$\alpha^{d}(c_{1})=y^{s}(c_{1})$ follows for some $1\leq s\leq r$. Now
$\gcd(e^{d}-1,r)=1$ implies that $(e^{d}-1)u+s=rv$ for some $u,v\in\mathbb{Z}$
and $y^{u}(c_{1})$ is a fixed point of $\alpha^{d}$. Indeed, the application
of Lemma 2.4 gives that%
\[
\alpha^{d}(y^{u}(c_{1}))=(\alpha^{d}\circ y^{u})(c_{1})=(y^{e^{d}u}\circ
\alpha^{d})(c_{1})=y^{e^{d}u}(\alpha^{d}(c_{1}))=y^{e^{d}u}(y^{s}(c_{1}))=
\]%
\[
y^{e^{d}u+s}(c_{1})=y^{u+rv}(c_{1})=y^{u}(c_{1}).
\]
Thus, $(y^{u}(c_{1}),\alpha(y^{u}(c_{1})),\ldots,\alpha^{d-1}(y^{u}(c_{1})))$
is a $d$-element cycle of $\alpha$. $\square$

\newpage

\noindent3.\thinspace NON-COMMUTING SOLUTIONS\thinspace OF\thinspace
$\alpha\circ y\circ\alpha^{-1}=y^{e}$

\bigskip

\noindent If $y\in\mathrm{S}_{n}$ is a solution of $\alpha\circ y\circ
\alpha^{-1}=y^{e}$, then the conditions $\alpha\circ y=y\circ\alpha$ and
$y^{e-1}=1$ are obviously equivalent to each other. A solution of
$(e\ast\varepsilon)$\ is called \textit{commuting} if one of these two
equivalent conditions holds. On the other hand, if $y\in S_{n}$\ is a
permutation with $y^{e-1}=1$ and $\alpha\circ y=y\circ\alpha$, then $y$\ is a
solution of $(e\ast\alpha)$: we have $y^{e}=y$ and $\alpha\circ y\circ
\alpha^{-1}=y=y^{e}$.

\bigskip

\noindent\textbf{3.1. Theorem.}\textit{ For any integer }$r\geq2$\textit{
dividing }$n\geq3$\textit{ and }$e^{\frac{n}{r}}-1$\textit{ and any }%
$n$\textit{-cycle }$\alpha\in\mathrm{S}_{n}$\textit{, the equation }%
$\alpha\circ y\circ\alpha^{-1}=y^{e}$\textit{ has a solution }$y\in
\mathrm{S}_{n}$\textit{, which is a product of }$n/r$\textit{\ disjoint }%
$r$\textit{-cycles. If the additional condition that }$r$\textit{ is not a
divisor of }$e-1$\textit{ holds, then }$y^{e-1}\neq1$\textit{ and }%
$y$\textit{\ is a non-commuting solution of }$(e\ast\alpha)$\textit{.}

\bigskip

\noindent\textbf{Proof.} Take $q=n/r$ and define an $n$-element set $P_{n}%
$\ of ordered pairs and a function $\varepsilon:P_{n}\longrightarrow P_{n}%
$\ as follows:%
\[
P_{n}=\{(i,j)\mid1\leq i\leq q\text{ and }1\leq j\leq r\},
\]%
\[
\varepsilon(i,j)=\left\{
\begin{array}
[c]{c}%
(i+1,ej)\text{ if }1\leq i\leq q-1\text{ and }1\leq j\leq r\\
(1,ej+1)\text{ if }i=q\text{ and }1\leq j\leq r\text{ \ \ \ \ \ \ \ \ \ \ \ }%
\end{array}
\right.  ,
\]
where $ej$ and $ej+1$ are taken in $\{1,2,\ldots,r\}$ modulo $r$. It is
straightforward to check that $\varepsilon$ is injective (hence a permutation
of $P_{n}$). Clearly, the definition of $\varepsilon$ immediately gives that
the length of a cycle in $\varepsilon$ is a multiple of $q$. If $1\leq k\leq
q-1$ and $1\leq j\leq r$, then $\varepsilon^{k}(1,j)=(1+k,e^{k}j)$ and $r\mid
e^{q}-1$ ensures that $\varepsilon^{q}(1,j)=(1,e^{q}j+1)=(1,j+1)$. Thus,%
\[
\varepsilon^{q}(1,j)=(1,j+1),\varepsilon^{2q}(1,j)=(1,j+2),\ldots
,\varepsilon^{(r-1)q}(1,j)=(1,j+r-1)
\]
are distinct elements and $\varepsilon^{rq}(1,j)=(1,j+r)=(1,j)$ (in $P_{n}$).
It follows that $\varepsilon$ has exactly one cycle of length $n=rq$. Now
define a permutation $y:P_{n}\longrightarrow P_{n}$\ as follows:%
\[
y(i,j)=(i,j+1),
\]
where $1\leq i\leq q$, $1\leq j\leq r$ and $j+1$ is taken in $\{1,2,\ldots
,r\}$ modulo $r$. The number of cycles in $y$ is exactly $q$ and each cycle of
$y$ is of the form%
\[
\left(  (i,1),(i,2),\ldots,(i,r)\right)  .
\]
An easy calculation shows that $\varepsilon\circ y=y^{e}\circ\varepsilon$. If
$1\leq i\leq q-1$ and $1\leq j\leq r$, then%
\[
\varepsilon(y(i,j))=\varepsilon((i,j+1))=(i+1,ej+e)
\]
and%
\[
y^{e}(\varepsilon(i,j))=y^{e}((i+1,ej))=(i+1,ej+e).
\]
If $i=q$ and $1\leq j\leq r$, then%
\[
\varepsilon(y(q,j))=\varepsilon((q,j+1))=(1,ej+e+1)
\]
and%
\[
y^{e}(\varepsilon(q,j))=y^{e}((1,ej+1))=(1,ej+1+e).
\]
Thus, $y$ is a required solution of $(e\ast\varepsilon)$. If $r$ is not a
divisor of $e-1$, then $y^{e-1}\neq1$. Finally, we deduce, that $(e\ast
\alpha)$ has similar solutions for any permutation $\alpha\in\mathrm{S}_{n}$
with $\mathrm{type}(\alpha)=\mathrm{type}(\varepsilon)$. $\square$

\bigskip

\noindent\textbf{3.2. Corollary.}\textit{ If }$\alpha\in\mathrm{S}_{n}%
$\textit{ is an }$n$\textit{-cycle and }$p\geq2$\textit{ is a prime divisor of
}$d=\gcd(n,e^{n}-1)\neq1$\textit{, then }$(e\ast\alpha)$\textit{ has a
solution }$y\in\mathrm{S}_{n}$\textit{, wich is a product of }$n/p$%
\textit{\ disjoint }$p$\textit{-cycles (and }$y^{d}=1$\textit{). If the
additional condition }$d\nmid(e-1)^{d}$\textit{ holds, then }$p$\textit{ can
be chosen as a non-divisor of }$e-1$\textit{. In this case }$y^{e-1}\neq
1$\textit{ and }$y$\textit{\ is a non-commuting solution of }$(e\ast\alpha)$.

\bigskip

\noindent\textbf{Proof.} Now $n=n_{1}p$ and%
\[
e^{\frac{n}{p}}-1=e^{n_{1}}-1=(e^{n_{1}p}-1)-e^{n_{1}}(e^{n_{1}(p-1)}%
-1)=(e^{n}-1)-e^{n_{1}}(e^{n_{1}(p-1)}-1)
\]
is divisible by $p$ (i.e. $p\mid e^{\frac{n}{p}}-1$) by Fermat's divisibility
$p\mid e^{p-1}-1$. The application of Theorem 3.1 gives the existence of a
solution $y\in\mathrm{S}_{n}$ of $(e\ast\alpha)$ such that $y$ is a product of
$n/p$\ disjoint $p$-cycles. Clearly, $p\mid d$ implies that $y^{d}=1$. In the
additional case $p\nmid e-1$ implies that $y^{e-1}\neq1$. $\square$

\bigskip

\noindent\textbf{3.3. Corollary.}\textit{ If }$\alpha\in\mathrm{S}_{n}%
$\textit{ is a permutation of }$\mathrm{type}(\alpha)=\left\langle g_{1}%
,g_{2},\ldots,g_{n}\right\rangle $\textit{ such that }$g_{a}\neq0$\textit{ and
}$p\geq2$\textit{ is a prime divisor of }$d=\gcd(a,e^{a}-1)\neq1$\textit{ for
some }$2\leq a\leq n$\textit{, then }$(e\ast\alpha)$\textit{ has a solution
}$y\in\mathrm{S}_{n}$\textit{ wich is a product of }$a/p$\textit{\ disjoint
}$p$\textit{-cycles (and }$y^{d}=1$\textit{). If the additional condition
}$d\nmid(e-1)^{d}$\textit{ holds, then }$p$\textit{ can be chosen as a
non-divisor of }$e-1$\textit{.} \textit{In this case }$y^{e-1}\neq1$\textit{
and }$y$\textit{\ is a non-commuting solution of }$(e\ast\alpha)$\textit{.}

\bigskip

\noindent\textbf{Proof.} If $(i,\alpha(i),\ldots,\alpha^{a-1}(i))$ is an
$a$-element cycle of $\alpha$, then $\{1,2,\ldots,n\}=H_{1}\cup H_{2}$, where
$H_{1}=\{i,\alpha(i),\ldots,\alpha^{a-1}(i)\}$ and $H_{2}=\{1,2,\ldots
,n\}\smallsetminus H_{1}$ are disjoint fixed sets of $\alpha$. Now Corollary
3.2 ensures the existence of a solution $y_{1}\in S_{H_{1}}$ of $(e\ast
(\alpha\upharpoonright H_{1}))$ (here $y_{1}:H_{1}\longrightarrow H_{1}$ is a
permutation and $(\alpha\upharpoonright H_{1}):H_{1}\longrightarrow H_{1}$ is
the restriction of $\alpha$\ to $H_{1}$) wich is a product of $a/p$\ disjoint
$p$-cycles. If $y_{2}=1$ is the identity permutation on $H_{2}$, then
$y=y_{1}\sqcup y_{2}$ is also a product of $a/p$\ disjoint $p$-cycles and
$(y_{1}\sqcup y_{2})^{d}=1$. The application of Lemma 2.3 gives that
$y\in\mathrm{S}_{n}$ is a solution of $(e\ast\alpha)$. In the additional case
Corollary 3.2 gives that $y_{1}$ is non-commuting, whence the non-commuting
property of $y_{1}\sqcup y_{2}$ follows. $\square$

\bigskip

\noindent\textbf{3.4. Proposition.}\textit{ Let }$\alpha=(1,2,\ldots
,n)\in\mathrm{S}_{n}$\textit{ be a cyclic permutation and }$r\geq2$\textit{ be
a divisor of }$n$\textit{ such that }$\gcd(\frac{n}{r},e^{n}-1)=1$\textit{. If
}$1\neq y\in\mathrm{S}_{n}$\textit{ is a solution of }$(e\ast\alpha)$\textit{
with }$\mathrm{type}(y)=\left\langle t_{1},t_{2},\ldots,t_{n}\right\rangle
$\textit{ and }$t_{s}\neq0$\textit{ for some }$2\leq s\leq n$\textit{, then
}$s$\textit{ is a divisor of }$r$\textit{ (i.e. any cycle length of }%
$y$\textit{\ is a divisor of }$r$\textit{), }$t_{s}$\textit{ is a multiple of
}$\frac{n}{r}$\textit{\ and }$y^{r}=1$\textit{. If }$r=p$\textit{ is a prime
number, then }$s=p$\textit{, }$t_{p}=\frac{n}{p}$\textit{ and }$t_{k}%
=0$\textit{ for each }$1\leq k\leq n$\textit{, }$k\neq p$\textit{, thus }%
$y$\textit{\ is a product of }$n/p$\textit{\ disjoint }$p$\textit{-cycles.}

\bigskip

\noindent\textbf{Proof.} Since $s\geq2$ is a cycle of length of $y$, Lemma 2.5
ensures that $s\mid e^{n}-1$, where $n=\mathrm{ord}(\alpha)$. Let $C_{i}%
^{(s)}\subseteq\{1,2,\ldots,n\}$, $1\leq i\leq t_{s}$\ be the pairwise
disjoint base sets of the $s$-element cycles of $y$. The type of $\alpha$\ and
Lemma 2.6\ ensure that $t_{s}s\in F_{1}(\alpha)=\{0,n\}=\{0,r\cdot\frac{n}%
{r}\}$. in view of $t_{s}s=r\cdot\frac{n}{r}$, $s\mid e^{n}-1$ and $\gcd
(\frac{n}{r},e^{n}-1)=1$, we deduce that $s\mid r$ and $\frac{n}{r}\mid t_{s}%
$. If $r=p$\ is prime, then $s=p$ and $t_{p}=\frac{n}{p}$\ follow from the
divisibility $s\mid p$ and $t_{s}s=p\cdot\frac{n}{p}$. $\square$

\bigskip

\noindent\textbf{3.5. Corollary.}\textit{ If }$\alpha=(1,2,\ldots
,n)\in\mathrm{S}_{n}$\textit{ is a cyclic permutation and }$p\geq2$\textit{ is
a prime divisor of }$n$\textit{ such that }$\gcd(\frac{n}{p},e^{n}%
-1)=1$\textit{ and }$p\nmid e-1$,\textit{ then }$y=1$\textit{ is the only
commuting solution of }$(e\ast\alpha)$\textit{.}

\bigskip

\noindent\textbf{Proof.} For any solution $1\neq y\in\mathrm{S}_{n}$\ of
$(e\ast\alpha)$\ Proposition 3.4 gives that $y$\ is a product of
$n/p$\ disjoint $p$-cycles. Now $y^{e-1}\neq1$ is a consequence of $p\nmid
e-1$. $\square$

\bigskip

\noindent\textbf{3.6. Theorem.}\textit{ Let }$\alpha=(1,2,\ldots
,n)\in\mathrm{S}_{n}$\textit{ be a cyclic permutation, }$p\geq2$\textit{ be a
prime divisor of }$n$\textit{ such that }$p\mid e^{\frac{n}{p}}-1$\textit{ and
}$\gcd(\frac{n}{p},e^{n}-1)=1$\textit{. Now the existence of a non-trivial
solution }$1\neq y\in\mathrm{S}_{n}$\textit{\ of }$(e\ast\alpha)$\textit{
follows from Theorem 3.1. If }$y\neq1$\textit{ is an arbitrary solution of
}$(e\ast\alpha)$\textit{, then the powers }$y,y^{2},\ldots,y^{p-1},y^{p}%
=1$\textit{ provide the complete list of solutions of }$(e\ast\alpha
)$\textit{.}

\bigskip

\noindent\textbf{Proof.} Let $1\neq y\in S_{n}$ be an arbitrary solution of
$(e\ast\alpha)$ with $\mathrm{type}(y)=\left\langle t_{1},t_{2},\ldots
,t_{n}\right\rangle $. Proposition 3.4 ensures that $y$\ is a product of
$q=n/p$\ disjoint $p$-cycles and the powers $y,y^{2},\ldots,y^{p-1},y^{p}=1$
are distinct solutions of $(e\ast\alpha)$.

\noindent If $1\leq d\leq t_{p}=q$ is the length of a cycle $(i_{1}%
,i_{2},\ldots,i_{d})$ of the induced permutation $\alpha^{(p)}$, then Lemma
2.6 gives $dp\in F_{d}(\alpha)\subseteq F_{1}(\alpha)=\{0,p\cdot q\}$ and
$d=q$. It follows that $\alpha^{(p)}$ has only one cycle of length $d=q$\ and
we can assume that $\alpha^{(p)}=(1,2,\ldots,q)$ (hence $\alpha(C_{j}%
^{(p)})=C_{j+1}^{(p)}$ for $1\leq j\leq q-1$ and $\alpha(C_{q}^{(p)}%
)=C_{1}^{(p)}$). Now fix an element $a\in C_{1}^{(p)}$. The above cyclic
property of $\alpha^{(p)}$\ and $1\leq q=n/p\leq n-1$ ensure that $a\neq
\alpha^{q}(a)\in C_{1}^{(p)}$, whence $\alpha^{q}(a)=y^{s}(a)$ follows for
some $1\leq s\leq p-1$ ($s$ depends on the choice of $a$). We claim that
$y(\alpha^{i}(a))=\alpha^{\overline{s}(i)q+i}(a)$ holds for all integers
$0\leq i\leq pq=n$, where $1\leq\overline{s}(i)\leq p-1$ denotes the
multiplicative inverse (reciprocal) of $e^{i}s$ in the prime field
$\mathbb{Z}_{p}$ (notice that $p\mid e^{i}s$ would contradict to $1\leq s\leq
p-1$ and $p\mid e^{\frac{n}{p}}-1$).

\noindent First we use $\alpha^{i}\circ y^{s}=y^{e^{i}s}\circ\alpha^{i}$ (in
Lemma 2.4) to get%
\[
\alpha^{q}(\alpha^{i}(a))=\alpha^{i}(\alpha^{q}(a))=\alpha^{i}(y^{s}%
(a))=y^{e^{i}s}(\alpha^{i}(a))
\]
and then we prove by induction, that%
\[
\alpha^{kq}(\alpha^{i}(a))=y^{ke^{i}s}(\alpha^{i}(a))
\]
holds for all integers $k\geq0$. Lemma 2.4 gives that $\alpha^{kq}\circ
y^{e^{i}s}=y^{e^{kq}e^{i}s}\circ\alpha^{kq}$ and $e^{kq}e^{i}s-e^{i}%
s=(e^{kq}-1)e^{i}s$ is divisible by $e^{q}-1$ and hence by $p$. Now
$e^{kq}e^{i}s-e^{i}s=pv$ and the validity of the above equality for $k+1$ can
be derived as%
\[
\alpha^{(k+1)q}(\alpha^{i}(a))=\alpha^{kq}(\alpha^{q}(\alpha^{i}%
(a)))=\alpha^{kq}(y^{e^{i}s}(\alpha^{i}(a)))=
\]%
\[
y^{e^{kq}e^{i}s}(\alpha^{kq}(\alpha^{i}(a)))=y^{e^{kq}e^{i}s}(y^{ke^{i}%
s}(\alpha^{i}(a)))=y^{e^{kq}e^{i}s-e^{i}s}(y^{(k+1)e^{i}s}(\alpha^{i}(a)))=
\]%
\[
y^{pv}(y^{(k+1)e^{i}s}(\alpha^{i}(a)))=y^{(k+1)e^{i}s}(\alpha^{i}(a)).
\]
In view of $\overline{s}(i)e^{i}s=1+pu$, the substitution $k=\overline{s}(i)$
into $\alpha^{kq}(\alpha^{i}(a))=y^{ke^{i}s}(\alpha^{i}(a))$ gives our claim%
\[
y(\alpha^{i}(a))=y^{\overline{s}(i)e^{i}s}(\alpha^{i}(a))=\alpha^{\overline
{s}(i)q}(\alpha^{i}(a))=\alpha^{\overline{s}(i)q+i}(a).
\]
Since $\{\alpha^{i}(a)\mid0\leq i\leq n-1\}=\{1,2,\ldots,n\}$, the permutation
$y$ is completely determined by $s$. The fact that $1\leq s\leq p-1$ can be
chosen in $p-1$ different ways implies that the number of non-trivial
solutions of $(e\ast\alpha)$ is at most $p-1$. It follows that $y,y^{2}%
,\ldots,y^{p-1},y^{p}=1$ is the complete list of solutions. $\square$

\bigskip

\noindent\textbf{3.7. Examples.} The following integers satisfy the conditions
in 3.1, 3.5 and 3.6. If $e=2$ and $n=20,21$, then we have $5\mid2^{\frac
{20}{5}}-1$, $\gcd(\frac{20}{5},2^{20}-1)=1$, $5\nmid2-1$\ and $7\mid
2^{\frac{21}{7}}-1$, $\gcd(\frac{21}{7},2^{21}-1)=1$, $7\nmid2-1$. If $e=-2$
and $n=55$, then we have $11\mid(-2)^{\frac{55}{11}}-1$, $\gcd(\frac{55}%
{11},2^{55}-1)=1$, $11\nmid(-2)-1$\textit{.}

\newpage

\noindent4.\thinspace COMMUTING\ SOLUTIONS\thinspace OF\thinspace$\alpha\circ
y\circ\alpha^{-1}=y^{e}$

\bigskip

\noindent\textbf{4.1. Theorem.}\textit{ Let }$w=\mathrm{ord}(\alpha
)=\operatorname{lcm}\{a\mid1\leq a\leq n,g_{a}\neq0\}$\textit{ be the order of
the permutation }$\alpha\in\mathrm{S}_{n}$\textit{ of }$\mathrm{type}%
(\alpha)=\left\langle g_{1},g_{2},\ldots,g_{n}\right\rangle $\textit{ with
}$g_{1}=0$\textit{. Assume that for any choice of the integers }$r\geq
2$\textit{ and }$d\geq2$\textit{ with }$\gcd(e-1,r)=1$\textit{, }$r\mid
e^{w}-1$\textit{, }$d\mid w$\textit{\ and }$dr\in F_{d}(\alpha)$\textit{ we
have }$g_{d}=0$\textit{ and }$\gcd(e^{d}-1,r)=1$\textit{. If }$y$\textit{ is a
solution of }$(e\ast\alpha)$\textit{ such that }$\gcd(e-1,s)=1$\textit{ for
any cycle length }$s\geq2$\textit{ of }$y$\textit{, then }$y=1$\textit{.}

\bigskip

\noindent\textbf{Proof.} Let $y\in S_{n}$ be a solution of $(e\ast\alpha)$
with $\mathrm{type}(y)=\left\langle t_{1},t_{2},\ldots,t_{n}\right\rangle $.
Assume that $y\neq1$, then $y$ has at least one cycle of length $r\geq2$. Now
we have $\gcd(e-1,r)=1$ and Lemma 2.5 ensures that $r\mid e^{w}-1$. Let
$C_{i}^{(r)}\subseteq\{1,2,\ldots,n\}$, $1\leq i\leq t_{r}$\ be the pairwise
disjoint base sets of the $r$-element cycles of $y$. In view of Lemma 2.6 we
have $\alpha(C_{i}^{(r)})=C_{\gamma(i)}^{(r)}$, $1\leq i\leq t_{r}$, where
$\gamma=\alpha^{(r)}$ is the induced permutation of the indices $\{1,2,\ldots
,t_{r}\}$. Consider an arbitrary cycle $(i_{1},i_{2},\ldots,i_{d}%
)$\ of$\ \alpha^{(r)}$, then $1\leq d\leq t_{r}$ and Lemma 2.6 gives $d\mid w$
and $dr\in F_{d}(\alpha)$. If $d\geq2$, then we have $g_{d}=0$ and $\gcd
(e^{d}-1,r)=1$. If $d=1$, then we also have $\gcd(e^{d}-1,r)=1$. Thus, the
application of Lemma 2.7 gives that $\alpha$\ also has a cycle of length $d$,
in contradiction with $g_{d}=0$. $\square$

\bigskip

\noindent\textbf{4.2. Corollary.}\textit{ Let }$\alpha=(1,2,\ldots
,a)\circ(a+1,a+2,\ldots,a+b)\in\mathrm{S}_{n}$\textit{ be a product of two
cyclic permutations, where }$2\leq a<b\leq n=a+b$\textit{. If }$a$\textit{ is
not a divisor of }$b$\textit{ and }$\gcd(u,e^{u}-1)=1$\textit{ for any choice
of }$u\in\{a,b,a+b\}$\textit{, then }$(e\ast\alpha)$\textit{ has only the
trivial solution }$y=1$\textit{.}

\bigskip

\noindent\textbf{Proof}. In order to use Theorem 4.1 take the integers
$r\geq2$ and $d\geq2$ such that

\noindent$\gcd(e-1,r)=1$, $r\mid e^{w}-1$, $d\mid w$ and%
\[
u=dr\in F_{d}(\alpha)\subseteq F_{1}(\alpha)=\{0,a,b,a+b\},
\]
where $w=\mathrm{ord}(\alpha)=\operatorname{lcm}(a,b)$. Clearly,
$d\notin\{a,b\}$ is a consequence of the fact that $a$ is not a divisor of
$b$. It follows that $g_{1}=g_{d}=0$ in $\mathrm{type}(\alpha)$ and
$\gcd(e^{d}-1,r)=1$ is a consequence of $\gcd(e^{dr}-1,dr)=\gcd(e^{u}-1,u)=1$.

\noindent If $y$ is a solution of $(e\ast\alpha)$ with $\mathrm{type}%
(y)=\left\langle t_{1},t_{2},\ldots,t_{n}\right\rangle $ and $s\geq2$ is a
cycle length of $y$, then $t_{s}s\in F_{1}(\alpha)=\{0,a,b,a+b\}$ by Lemma
2.6. In view of $\gcd(e-1,a)=\gcd(e-1,b)=\gcd(e-1,a+b)=1$, we obtain that
$\gcd(e-1,s)=1$. Thus Theorem 4.1 ensures that $y=1$\textit{. }$\square$

\bigskip

\noindent\textbf{4.3. Corollary.}\textit{ Let }$\alpha=(1,2,\ldots
,a)\circ(a+1,a+2,\ldots,a+b)\in\mathrm{S}_{n}$\textit{ be a product of two
cyclic permutations, where }$2\leq a\leq b\leq n=a+b$\textit{. If }$a$\textit{
is a divisor of }$b$\textit{, }$\gcd(b,e^{b}-1)=\gcd(a+b,e-1)=1$\textit{\ and
}$\frac{a+b}{a}$\textit{ is not a divisor of }$e^{b}-1$\textit{, then }%
$(e\ast\alpha)$\textit{ has only the trivial solution }$y=1$\textit{.}

\bigskip

\noindent\textbf{Proof}. In order to use Theorem 4.1 take the integers
$r\geq2$ and $d\geq2$ such that

\noindent$\gcd(e-1,r)=1$, $r\mid e^{w}-1$, $d\mid w$ and%
\[
u=dr\in F_{d}(\alpha)\subseteq F_{1}(\alpha)=\{0,a,b,a+b\},
\]
where $w=\mathrm{ord}(\alpha)=\operatorname{lcm}(a,b)$. The fact that $a$\ is
a divisor of $b$ implies that $w=b$. Since $\frac{a+b}{a}$ is not a divisor of
$e^{b}-1$, there are no integers $r\geq2$ and $d\geq2$ with the above
properties. Indeed, $dr\in\{a,b\}$ is in contradiction with $r\mid e^{b}-1$,
$r\mid b$ and $\gcd(e^{b}-1,b)=1$. If $dr=a+b$, then $d\mid b$ implies that
$d\mid a$ and $\frac{a+b}{a}\mid\frac{a+b}{d}$. Now $r=\frac{a+b}{d}$ is in
contradiction with $\frac{a+b}{a}\nmid e^{b}-1$. Thus, we do not have to check
$g_{d}=0$ and $\gcd(e^{d}-1,r)=1$.

\noindent If $y$ is a solution of $(e\ast\alpha)$ with $\mathrm{type}%
(y)=\left\langle t_{1},t_{2},\ldots,t_{n}\right\rangle $ and $s\geq2$ is a
cycle length of $y$, then $t_{s}s\in F_{1}(\alpha)=\{0,a,b,a+b\}$ by Lemma
2.6. In view of $\gcd(e-1,a)=\gcd(e-1,b)=\gcd(e-1,a+b)=1$, we obtain that
$\gcd(e-1,s)=1$. Thus Theorem 4.1 ensures that $y=1$\textit{. }$\square$

\bigskip

\noindent\textbf{4.4. Examples.} If $e=2$, then $a=10$ and $b=15$ satisfy the
conditions in Corollary 4.2. If $e=-2$, then $a=35$ and $b=77$ satisfy the
conditions in Corollary 4.2. If $e=2$, then $a=5$ and $b=25$ satisfy the
conditions in Corollary 4.3. If $e=-2$, then $a=2$ and $b=14$ satisfy the
conditions in Corollary 4.3.

\bigskip

\noindent For an integer $2\leq v$\ with $\gcd(v,e-1)=1$ let $q(e,v)$ denote
the smallest prime divisor of $e^{v}-1$ not dividing $e-1$ (if there is no
such prime, then take $q(e,v)=+\infty$). If $v$\ is odd, then $7\leq q(2,v)$.
Similarly, if $v$\ is not divisible by $2$\ and $3$, then $23\leq q(2,v)$.
Indeed, $2\nmid v$ and $3\nmid v$\ imply that $2,3,5,7,11,13,17,19\nmid
2^{v}-1$. Notice that $q(2,2)=2^{2}-1$, $q(2,3)=2^{3}-1$, $q(2,4)=3$,
$q(2,5)=2^{5}-1$, $q(2,6)=3$, $q(2,7)=2^{7}-1$, $q(2,8)=3$, $q(2,9)=7$,
$q(2,10)=3$ and $q(2,11)=23\neq2^{11}-1$ (clearly, $2^{11}-1$ is not a
Mersenne prime). We also note that $q(-2,2)=+\infty$, $q(-2,4)=5$,
$q(-2,5)=11$, $q(-2,7)=43$, $q(-2,8)=5$, $q(-2,10)=11$ and $q(-2,11)=683$.

\bigskip

\noindent\textbf{4.5. Theorem.}\textit{ Let }$w=\mathrm{ord}(\alpha
)=\operatorname{lcm}\{a\mid1\leq a\leq n,g_{a}\neq0\}$\textit{ be the order of
the permutation }$\alpha\in\mathrm{S}_{n}$\textit{ of type }$\mathrm{type}%
(\alpha)=\left\langle g_{1},g_{2},\ldots,g_{n}\right\rangle $\textit{ such
that }$g_{1}=0$\textit{ and }$\gcd(a,b)=1$\textit{ for all }$1\leq a<b\leq
n$\textit{ with }$g_{a}\neq0\neq g_{b}$\textit{ (two different cycle lengths
are relative primes). Assume that }$1\leq g_{a}\leq q(e,w)-1$\textit{ and
}$\gcd(a,e^{a}-1)=1$\textit{ hold for all }$1\leq a\leq n$\textit{ with
}$g_{a}\neq0$\textit{. If }$y\in\mathrm{S}_{n}$\textit{\ is a solution of
}$(e\ast\alpha)$\textit{, then }$y^{e-1}=1$\textit{ (i.e. }$y$\textit{ is
commuting).}

\bigskip

\noindent\textbf{Proof}. Let $y\in S_{n}$ be a solution of $(e\ast\alpha)$.
Any cycle $C$\ of $y^{(e-1)^{n}}$ can be obtained from a uniquely determined
cycle $D$\ of $y$ (notice that $C\subseteq D$). If $s=\left\vert D\right\vert
$ is the length of $D$, then $\left\vert C\right\vert =\frac{s}{d}$ with
$d=\gcd((e-1)^{n},s)$. Since $\gcd(e-1,\frac{s}{d})=1$, we can apply Theorem
4.1 on $y^{(e-1)^{n}}$. Let $r\geq2$ and $d\geq2$ be integers such that
$\gcd(e-1,r)=1$, $r\mid e^{w}-1$, $d\mid w$\textit{\ }and $dr\in F_{d}%
(\alpha)$. Our conditions on the type of $\alpha$\ ensure the existence of a
unique integer $k\geq1$ such that $1\leq dk\leq n$ and $g_{dk}\neq0$, whence
$F_{d}(\alpha)=F_{dk}(\alpha)=\{0,dk,2dk,\ldots,g_{dk}dk\}$ follows. Now we
have $dr=jdk$ as well as $r=jk$ for some $1\leq j\leq g_{dk}$. Since $r\mid
e^{w}-1$, we obtain that $j\mid e^{w}-1$. Clearly, $\gcd(e-1,j)=1$ and $2\leq
j\leq g_{dk}\leq q(e,w)-1$ would contradict to the definition of $q(e,w)$. It
follows that $j=1$ and $g_{dr}=g_{dk}\neq0$. Thus, $g_{d}=0$ and $\gcd
(e^{dr}-1,dr)=1$\ are also consequences of our conditions on the type of
$\alpha$, whence $\gcd(e^{d}-1,r)=1$\ follows. The application of Theorem 4.1
gives that $y^{(e-1)^{n}}=1$. We also have $y^{e^{w}-1}=1$ by Lemma 2.5. Since
$\gcd(a,e-1)=1$\ holds for all $1\leq a\leq n$ with $g_{a}\neq0$\ (our
assumption that $\gcd(a,e^{a}-1)=1$ is stronger) and $w=\operatorname{lcm}%
\{a\mid1\leq a\leq n,g_{a}\neq0\}$, we obtain that $\gcd(w,e-1)=1$. In view of%
\[
e^{w}-1=(e-1)(w+(e-1)+(e^{2}-1)+\cdots+(e^{w-1}-1)),
\]
$\gcd(e^{w}-1,(e-1)^{n})=e-1$ can be derived, whence $(e^{w}-1)u+(e-1)^{n}%
v=e-1$ follows for some $u,v\in\mathbb{Z}$. As a consequence we have%
\[
y^{e-1}=y^{(e^{w}-1)u+(e-1)^{n}v}=(y^{e^{w}-1})^{u}\circ(y^{(e-1)^{n}}%
)^{v}=1.
\]
Thus, $y$\ is a commuting solution of $(e\ast\alpha)$. $\square$

\bigskip

\noindent\textbf{4.6. Examples.} If we omit any one of the three conditions
$\gcd(a,b)=1$,

\noindent$1\leq g_{a}\leq q(e,w)-1$ and $\gcd(a,e^{a}-1)=1$\ from Theorem 4.5,
then we loose the validity of 4.5.

\noindent(i) The cycle lengths $2$ and $4$\ in $\alpha=(1,2)\circ
(3,4,5,6)\circ(7,8,9,10)\in\mathrm{S}_{10}$ are not relative primes, while the
other two conditions with $e=2$\ in 4.5 hold for $\alpha$ ($w=\mathrm{ord}%
(\alpha)=4$, $q(2,4)=3$ and $\gcd(2,2^{2}-1)=\gcd(4,2^{4}-1)=1$). Now it is
straightforward to check that $y=(1,10,4,6,8)\circ(2,3,7,9,5)\in
\mathrm{S}_{10}$ is a non-commuting solution of $(e\ast\alpha)$: $y^{e-1}%
\neq1$.

\noindent(ii) The only cycle length is $2$ in $\alpha=(1,2)\circ
(3,4)\circ(5,6)\in\mathrm{S}_{6}$ and $g_{2}=3\nleq q(2,2)-1=2$, while the
other two conditions with $e=2$\ in 4.5 hold for $\alpha$ ($w=\mathrm{ord}%
(\alpha)=2$ and $\gcd(2,2^{2}-1)=1$). Now it is straightforward to check that
$y=(1,3,5)\circ(2,6,4)\in\mathrm{S}_{6}$ is a non-commuting solution of
$(e\ast\alpha)$: $y^{e-1}\neq1$.

\noindent(iii) If $d=\gcd(a,e^{a}-1)\neq1$ for some $1\leq a\leq n$ with
$g_{a}\neq0$, then the weak additional condition $d\nmid(e-1)^{d}$ in
Corollary 3.3 ensures the existence of a non-commuting solution of
$(e\ast\alpha)$.

\bigskip

\noindent\textbf{4.7. Remarks.} Theorem 4.5 reveals that under certain mild
assumptions, the solutions of $(e\ast\alpha)$ are exactly the solutions of
$y^{e-1}=1$ in the centralizer of $\alpha$. If $a=p^{t}$ is a prime power and
$p\nmid e-1$, then it is straightforward to check that%
\[
e=e^{p^{0}}\equiv e^{p^{1}}\equiv\cdots\equiv e^{p^{t-1}}\equiv e^{p^{t}%
}\operatorname{mod}(p),
\]
whence $\gcd(a,e^{a}-1)=1$ follows. As $\gcd(w,e-1)$ appears in the above
proof of 4.5, we add the following simple observation. If $\alpha\in
\mathrm{S}_{n}$ is arbitrary and $\gcd(w,e-1)=d\neq1$, then $y=\alpha
^{\frac{w}{d}}\neq1$ is a non-trivial commuting solution of $(e\ast\alpha)$
such that $y^{d}=1$.

\bigskip

\noindent The replacement of $1\leq g_{a}\leq q(e,w)-1$ in Theorem 4.5 by the
stronger condition $g_{a}=1$ helps to give a complete answer when our equation
$(e\ast\alpha)$\ has only the trivial solution $y=1$. If $g_{1}=0$ and $0\leq
g_{a}\leq1$, then our next result shows that $\gcd(a,e^{a}-1)\neq1$\ in
Corollary 3.3 is essential for the existence of a non-trivial solution.

\bigskip

\noindent\textbf{4.8. Theorem.}\textit{ Let }$\alpha\in\mathrm{S}_{n}$\textit{
be a permutation such that }$\mathrm{type}(\alpha)=\left\langle g_{1}%
,g_{2},\ldots,g_{n}\right\rangle $\textit{, }$g_{1}=0$\textit{, }$0\leq
g_{a}\leq1$\textit{ and }$\gcd(a,b)=1$\textit{ for all }$1\leq a,b\leq
n$\textit{ with }$a\neq b$\textit{ and }$g_{a}=g_{b}=1$\textit{ (two different
cycle lengths are relative primes). Now the following conditions are
equivalent:}

\noindent\textit{(i) }$\gcd(a,e^{a}-1)=1$\textit{ holds for all }$1\leq a\leq
n$\textit{ with }$g_{a}=1$\textit{.}

\noindent\textit{(ii) The only solution }$y\in\mathrm{S}_{n}$\textit{\ of
}$(e\ast\alpha)$\textit{ is }$y=1$\textit{.}

\bigskip

\noindent\textbf{Proof}. (i)$\Longrightarrow$(ii): Now $\alpha=\alpha_{1}%
\circ\alpha_{2}\circ\cdots\circ\alpha_{k}$\ is a product of pairwise disjoint
cycles $\alpha_{1},\alpha_{2},\ldots,\alpha_{k}\in S_{n}$ such that%
\[
2\leq a_{1}=\mathrm{ord}(\alpha_{1})<a_{2}=\mathrm{ord}(\alpha_{2}%
)<\cdots<a_{k}=\mathrm{ord}(\alpha_{k})\leq n
\]
and $a_{1}+a_{2}+\cdots+a_{k}=n$. If $y\in\mathrm{S}_{n}$ is a solution of
$(e\ast\alpha)$, then Theorem 4.5 gives that $y$ is commuting: $y^{e-1}=1$ and
$\alpha\circ y=y\circ\alpha$. The type of $\alpha$\ (no two cycles of the same
length) ensures that
\[
\mathrm{Cen}(\alpha)=\{\alpha_{1}^{t_{1}}\circ\alpha_{2}^{t_{2}}\circ
\cdots\circ\alpha_{k}^{t_{k}}\mid0\leq t_{i}\leq a_{i}-1\text{ for each }1\leq
i\leq k\}.
\]
In view of $y\in\mathrm{Cen}(\alpha)$, we have $y=\alpha_{1}^{t_{1}}%
\circ\alpha_{2}^{t_{2}}\circ\cdots\circ\alpha_{k}^{t_{k}}$ for some integers
$0\leq t_{i}\leq a_{i}-1$, $1\leq i\leq k$. Clearly, $y^{e-1}=1$ implies that
the divisibility $a_{i}\mid(e-1)t_{i}$ holds for each $1\leq i\leq k$. Since
$\gcd(a_{i},e-1)=1$ is a consequence of $\gcd(a_{i},e^{a_{i}}-1)=1$, we deduce
that $t_{i}=0$ for each $1\leq i\leq k$. It follows that $y=1$.

\noindent(ii)$\Longrightarrow$(i): If $\gcd(a,e^{a}-1)\neq1$ holds for some
$1\leq a\leq n$ with $g_{a}=1$ and $p\geq2$ is a prime divisor of
$\gcd(a,e^{a}-1)$, then Corollary 3.3 gives the existence of a solution $1\neq
y\in S_{n}$\ of $(e\ast\alpha)$ which is a product of $a/p$\ disjoint
$p$-cycles. $\square$

\bigskip

\noindent\textbf{4.9. Remarks.} The following observations are closely related
to some of the conditions in Theorems 4.5 and 4.8. If $\alpha\in S_{n}$ is a
permutation such that $\mathrm{type}(\alpha)=\left\langle g_{1},g_{2}%
,\ldots,g_{n}\right\rangle $ and $g_{1}=0$, then we have the following
inequalities:%
\[
w=\mathrm{ord}(\alpha)\leq\underset{1\leq a\leq n}{\prod}a^{g_{a}}%
\leq\left\vert \mathrm{Cen}(\alpha)\right\vert .
\]
Since $w=\mathrm{ord}(\alpha)=\operatorname{lcm}\{a\mid1\leq a\leq n,g_{a}%
\neq0\}$, it is clear that $\mathrm{ord}(\alpha)=\underset{1\leq a\leq
n}{\prod}a^{g_{a}}$ holds if and only if $0\leq g_{a}\leq1$ and $\gcd(a,b)=1$
for all $1\leq a,b\leq n$ with $a\neq b$ and $g_{a}=g_{b}=1$. It is not
difficult to show that $\underset{1\leq a\leq n}{\prod}a^{g_{a}}=\left\vert
\mathrm{Cen}(\alpha)\right\vert $ holds if and only if $0\leq g_{a}\leq1$ for
all $1\leq a\leq n$.

\bigskip

\noindent\textbf{Acknowledgment}

\noindent The authors would like to sincerely thank the several insightful
comments and suggestions of the referee.

\bigskip

\noindent REFERENCES

\bigskip

\noindent\lbrack CHS] S. Chowla, I.N. Herstein, W.R. Scott:\textit{ The
solutions of }$x^{d}=1$\textit{ in symmetric groups, }Norske Vid. Selsk.
(Trindheim) 25, 29-31 (1952)

\noindent\lbrack E] A. Evangelidou: \textit{The Solution of Length Five
Equations Over Groups,} Communications in Algebra, 35:6, 1914-1948 (2007)

\noindent\lbrack FM] H.Finkelstein, K.I. Mandelberg:\textit{ On Solutions of
"Equations in Symmetric Groups"}, Journal of Combinatorial Theory, Series A,
25, 142-152 (1978)

\noindent\lbrack GR] M. Gerstenhaber, O. S. Rothaus: \textit{The solution of
sets of equations in groups}, Proc. Nat. Acad. Sci. U.S.A. 48, 1531-1533 (1962)

\noindent\lbrack IR] I. M. Isaacs, G. R. Robinson: \textit{On a Theorem of
Frobenius: Solutions of }$x^{n}=1$\textit{ in Finite Groups,} American
Mathematical Monthly, Vol. 99, No. 4, pp. 352-354 (Apr., 1992)

\noindent\lbrack I] I. Martin Isaacs:\textit{ Finite Group Theory,} American
Mathematical Society, 2008

\noindent\lbrack L] R. C. Lyndon: \textit{Equations in groups}, Boletim Da
Sociedade Brasileira de Matem\'{a}tica, 11(1), 79--102 (1980)

\noindent\lbrack MW] L. Moser, M. Wyman:\textit{ On solutions of }$x^{d}%
=1$\textit{ in symmetric groups,} Canad. J. Math. 7, 159-168 (1955)
\end{document}